\documentclass{article}
\usepackage{graphics, graphicx, latexsym, amssymb}

\begin{document}

\title{Stacked directed animals and multi-directed animals \\
defined without using heaps of pieces}
\author{Svjetlan Fereti\'{c} \footnote{e-mail: svjetlan.feretic@gradri.hr} \\ 
Faculty of Civil Engineering, University of Rijeka, \\ 
Viktora Cara Emina 5, 51000 Rijeka, Croatia}
\maketitle

\begin{abstract}
Stacked directed animals and multi-directed animals are two lattice models defined by Bousquet-M\'{e}lou and Rechnitzer in 2002. The original definitions of those models involve heaps of pieces, \textit{i.e.}, some geometric representation of partially commutative monoids. The object of this writing is to define stacked directed animals and multi-directed animals in such a way that heaps of pieces are not involved. Our alternative definitions are equivalent to Bousquet-M\'{e}lou and Rechnitzer's original ones.
\end{abstract}

\section{Introduction}

\textit{Directed animals} are a lattice model invented in 1982 by the physicists D. Dhar, M. K. Phani and M. Barma \cite{Dhar}. A \textit{heap of pieces} is an object that gives a convenient geometric representation of a partially commutative monoid. Heaps of pieces were invented in 1985 by the mathematician G. X. Viennot \cite{Montreal}. Viennot also observed that directed animals are in bijection with \textit{pyramids} of \textit{dimers} \cite{Bourbaki}. (A pyramid is a heap having only one minimal piece and a dimer is a horizontal, one unit long line segment connecting two lattice points.) \textit{Stacked directed animals} are a superset of directed animals. Stacked directed animals, as well as \textit{multi-directed animals} (a superset of stacked directed animals), were defined in 2002 by M. Bousquet-M\'{e}lou and A. Rechnitzer \cite{Rechnitzer}. In \cite{Rechnitzer}, the authors begin with defining \textit{stacked pyramids} and \textit{connected heaps}; connected heaps are a superset of stacked pyramids, which are, in turn, a superset of pyramids of dimers. Then they define an injection (say $\varphi$) from connected heaps into polyominoes; this injection is an extension of Viennot's bijection. Stacked directed animals are then defined to be the image of stacked pyramids under the injection $\varphi$, and multi-directed animals are defined to be the image of connected heaps under the injection $\varphi$. 

The purpose of this writing is to define stacked directed animals and multi-directed animals in a direct way, \textit{i.e.}, without using heaps of pieces. The original definition of directed animals does not involve heaps of pieces. In this sense, our alternative definitions homologize stacked directed animals and multi-directed animals with directed animals. Sure enough, these alternative definitions are equivalent to the definitions given by Bousquet-M\'{e}lou and Rechnitzer.

\section{Everything else}

There are three regular tilings of the Euclidean plane, namely the triangular tiling, the square tiling, and the hexagonal tiling. We adopt the convention that every square or hexagonal tile has two horizontal edges. In a regular tiling, a tile is often referred to as a \textit{cell}. A plane figure $P$ is a \textit{polyomino} if $P$ is a union of finitely many cells and the interior of $P$ is connected. See Figure 1. Observe that, if a union of \textit{hexagonal} cells is connected, then it possesses a connected interior as well.

\begin{figure}
\begin{center}
\includegraphics[width=55mm]{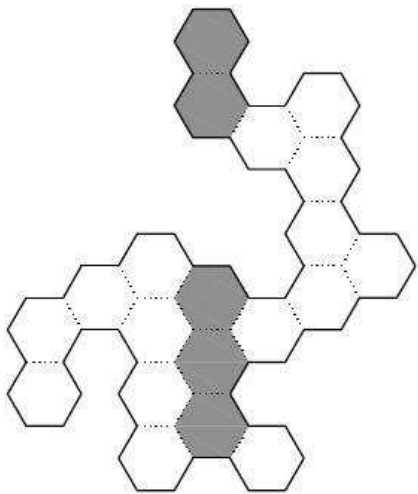}
\caption{A hexagonal-celled polyomino.}
\end{center}
\end{figure}

Let $P$ and $Q$ be two polyominoes. We consider $P$ and $Q$ to be equal if and only if there exists a translation $f$ such that $f(P)=Q$.

From now on, we concentrate on the hexagonal tiling. When we write ``a polyomino", we actually mean ``a hexagonal-celled polyomino".

Given a polyomino $P$, it is useful to partition the cells of $P$ according to their horizontal projection. Each block of that partition is a \textit{column} of $P$. Note that a column of a polyomino is not necessarily a connected set. An example of this is the highlighted column in Figure~1. On the other hand, it may happen that every column of a polyomino $P$ is a connected set. In this case, the polyomino $P$ is a \textit{column-convex polyomino}. See Figure 2.

\begin{figure}
\begin{center}
\includegraphics[width=55mm]{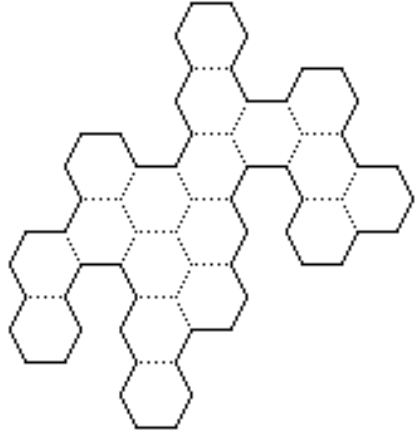}
\caption{A column-convex polyomino.}
\end{center}
\end{figure}

Let $c$ and $d$ be two cells on the hexagonal lattice, and let $c$ be either the upper left neighbour of $d$, or the upper neighbour of $d$, or the upper right neighbour of $d$. Then we say that $c$ is a \textit{generalized upper neighbour} of $d$.

\textit{Directed animals} can be defined by recursion on the number of cells. The one-celled polyomino is a directed animal. Let $B$ be a union of $n \geq 2$ cells. Then $B$ is a directed animal if and only if $B$ can be written as $B=A \cup c$, where $A$ is a directed animal with $n-1$ cells, and $c$ is a generalized upper neighbour of some cell of $A$. See Figure 3.

\begin{figure}
\begin{center}
\includegraphics[width=55mm]{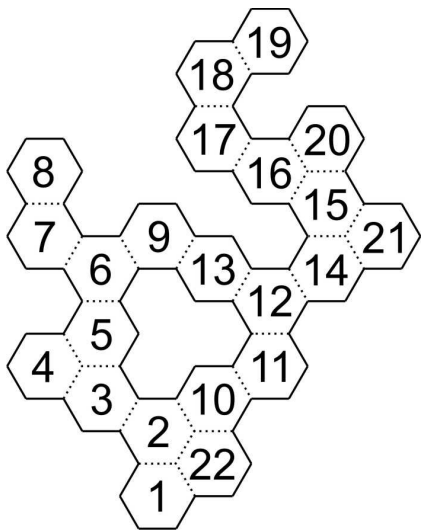}
\caption{A directed animal. For every $i \in \{2,\ 3,\ldots,\ 22\}$, the cells labeled $1,\ 2,\ldots,\ i-1$ form a directed animal, and the cell labeled $i$ is a generalized upper neighbour of some cell of the said directed animal.}
\end{center}
\end{figure}

Let $A$ be a directed animal. Among the cells of $A$, there exists exactly one cell which is not a generalized upper neighbour to any cell of $A$. That special cell is called the \textit{source cell} of $A$. 

We say that a cell $c$ \textit{dominates} over a cell $d$ if, for some $i \in \{0,\ 1,\ 2,\ldots \}$, $c$ lies $i$ units above a generalized upper neighbour of $d$. Suppose that $P$ is a polyomino, that $c$ is a cell of $P$, and that $c$ does not dominate over any cell of $P$. Then we say that $c$ is a \textit{source cell} of $P$. For directed animals, this definition of source cells is equivalent to the one given just above.

Let $P$ and $Q$ be two polyominoes that have no cells in common. We say that $P$ \textit{dominates} over $Q$ if there exists a cell of $P$ which dominates over a cell of $Q$, but there does not exist a cell of $Q$ which dominates over a cell of $P$.

A polyomino $P$ is a \textit{multi-directed animal} if there exist $k \in \mathbb{N}$ and directed animals $A_1,\ A_2,\ldots ,\ A_k$, no two of which have a cell in common, such that

\begin{itemize}
\item $P = \cup_{i=1}^k A_i$;
\item For $j \in \{2,\ 3,\ldots ,\ k\}$, the horizontal projection of the source cell of $A_{j-1}$ lies to the left of, and hence does not intersect with, the horizontal projection of the animal $A_{j}$;
\item However, for $j \in \{2,\ 3,\ldots ,\ k\}$, the union $\cup_{i=1}^{j-1} A_i$ dominates over the animal $A_j$;
\item For $j \in \{2,\ 3,\ldots ,\ k\}$, the union $\cup_{i=1}^{j-1} A_i$ and the animal $A_j$ have at least one edge in common.
\end{itemize}

In the definition of a \textit{stacked directed animal}, in addition to the four requirements stated above, there is also a fifth requirement, namely

\begin{itemize}
\item For $j \in \{2,\ 3,\ldots ,\ k\}$, $A_{j-1}$ itself dominates over $A_j$.
\end{itemize}

Thus, stacked directed animals are a subset of multi-directed animals. See Figures 4 and 5. The multi-directed animal of Figure 4 is not a stacked directed animal because $A_3$ does not dominate over $A_4$.

\begin{figure}
\begin{center}
\includegraphics[width=137.5mm]{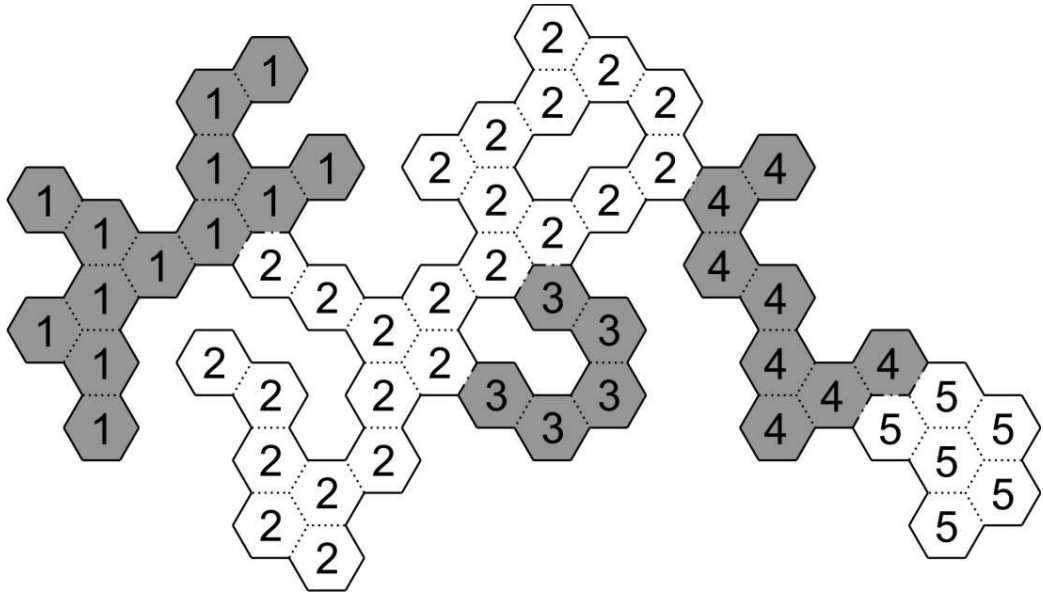}
\caption{A multi-directed animal. For $i \in \{1,\ldots ,\ 5\}$, the cells of $A_i$ are labeled by $i$.}
\end{center}
\end{figure}

\begin{figure}
\begin{center}
\includegraphics[width=5.5in]{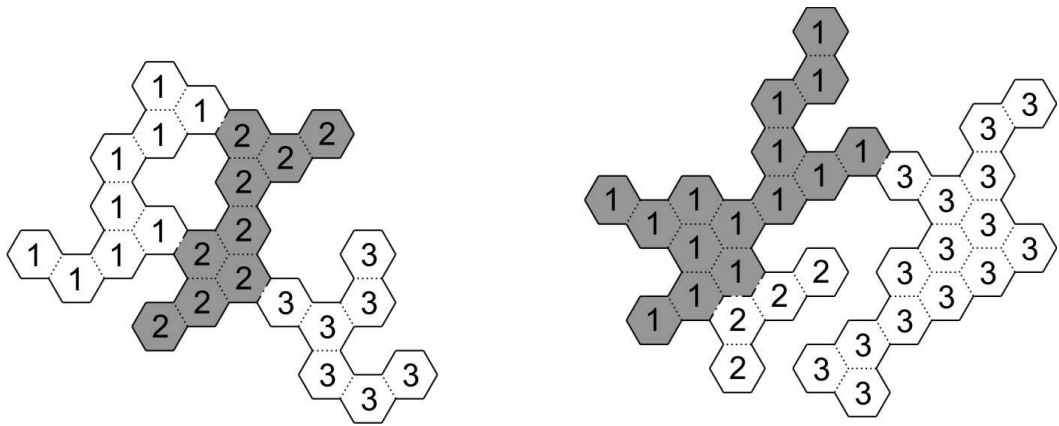}
\caption{Two stacked directed animals. For $i \in \{1,\ 2,\ 3\}$, the cells of $A_i$ are labeled by $i$.}
\end{center}
\end{figure}

For every multi-directed animal $P$, the positive integer $k$ and the directed animals $A_1,\ A_2,\ldots ,\ A_k$ (appearing in the definition of multi-directed animals) are unique. Namely, $k$ is the number of source cells of $P$. Let $c_1$ stand for the leftmost source cell of $P$, $c_2$ for the second-leftmost source cell of $P$, and so on. Then $A_1$ is the greatest instance of a directed animal which is contained in $P$ and has $c_1$ as the source cell; $A_2$ is the greatest instance of a directed animal which is contained in $P \setminus A_1$ and has $c_2$ as the source cell; $A_3$ is the greatest instance of a directed animal which is contained in $P \setminus (A_1 \cup A_2)$ and has $c_3$ as the source cell, \textit{etc.}

Every column-convex polyomino is a multi-directed animal; this can be proved by induction on the number of columns. By way of example, in Figure 6, the column-convex polyomino of Figure 2 is decomposed into three directed animals, and the decomposition has all the properties required in the definition of multi-directed animals. However, observe that the column-convex polyomino of Figure 2 \textit{is not} a stacked directed animal. 

\begin{figure}
\begin{center}
\includegraphics[width=55mm]{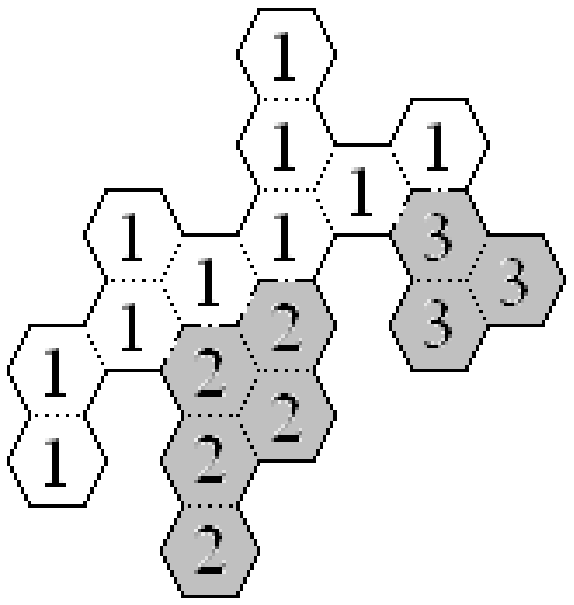}
\caption{The column-convex polyomino of Figure 2 decomposed into directed animals.}
\end{center}
\end{figure}

Likewise, there exist level one cheesy polyominoes\footnote{Cheesy polyominoes are an object defined by this author \cite{ja}.} that are not multi-directed animals. For example, if the level one cheesy polyomino of Figure 7 were a multi-directed animal, then it would be a union of two directed animals (because there are two source cells). The first directed animal would consist of the cells labeled by $1$, and the second directed animal of the cells labeled by $2$. The cell labeled by $e$ would not belong to any of the two directed animals, and hence would not belong to the cheesy polyomino. Since the cell $e$ does belong to the cheesy polyomino, we have a contradiction.

\begin{figure}
\begin{center}
\includegraphics[width=40mm]{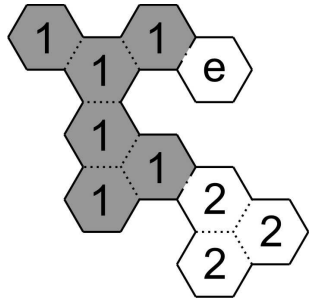}
\caption{A level one cheesy polyomino which is not a multi-directed animal.}
\end{center}
\end{figure}

On the other hand, there exist a lot of multi-directed animals and stacked directed animals which are not cheesy polyominoes. For example, none of the objects shown in Figures 4 and 5 is a cheesy polyomino.

\end{document}